\documentclass[a4paper,10pt]{article}
\usepackage[utf8]{inputenc}
\usepackage[affil-it]{authblk}

\usepackage{graphicx, caption, subcaption}
\usepackage{tikz}
\usetikzlibrary{matrix,arrows, positioning}
\usetikzlibrary[decorations.pathmorphing,mindmap]
\usetikzlibrary{mindmap,matrix,arrows,decorations.pathmorphing,patterns,fadings}

\usepackage{subfig, pgf}
\usepackage{amsfonts, amssymb, amsmath, amsthm,url}
\usepackage[normalem]{ulem}

\newcommand{\lra}{\longrightarrow}
\newcommand{\sm}{\setminus}
\newcommand{\ol}[1]{\overline{#1}}
\newcommand{\gal}{\textrm{Gal}(\ol{\QQ}/\QQ)}
\newcommand{\PP}{{\mathbf P}}

\newcommand{\CC}{{\mathbf C}}
\newcommand{\ZZ}{{\mathbf Z}}
\newcommand{\QQ}{{\mathbf Q}}

\newcommand{\lattes}{Latt\`{e}s}

\newtheorem{theorem}{Theorem}[section] % 1st argument is your name for it
\newtheorem{lemma}[theorem]{Lemma}     % 2nd argument is what is printed

\title{Belyi Latt\`{e}s maps}
\author{Ayberk Zeytin
	\thanks{Electronic address: \texttt{ayberkz@gmail.com}}}
\affil{Department of Mathematics, Galatasaray University \\ \c{C}{\i}ra\u{g}an Cad. No. 36, 34357 Be\c{s}ikta\c{s} \\\.{I}stanbul Turkey}

\begin{document}

\maketitle

\begin{abstract}
	In this work, we determine all {\lattes} maps which are Belyi morphisms. It turns out that in the generic case, i.e. when the automorphism group is $\ZZ/2\ZZ$, the corresponding family of {\lattes} maps are Belyi morphisms if and only if the isogeny is multiplication by two. This family form a continuous family of Belyi maps. Elliptic curves with complex multiplication also determine a family over $\ZZ$ of Belyi morphisms. We give the explicit formulas for the first few Belyi morphisms when the curve has complex multiplication by $3^{{rd}}$ root of unity.
\end{abstract}

\section{Introduction}
\label{sec:introduction}

An algebraic curve, $X$, admitting a model whose defining equations have algebraic coefficients is referred to as an \emph{arithmetic curve}. A celebrated theorem of Belyi, \cite{belyi}, states that such a curve admits a \emph{Belyi morphism}, $\phi$; that is, a meromorphic function ramified at most over $3$ points, which are assumed in literature to be $0$, $1$ and $\infty$. The pair $(X,\phi)$ is called a \emph{Belyi pair}. 

The absolute Galois group, $\gal$, acts on the set of all arithmetic curves equipped with Belyi morphisms or on Belyi pairs $(X,\phi)$ (see \cite[Section~2.4.1.2]{lando-zvonkine-lowdimtop} for a precise description). When one restricts to genus $0$ this reduces to understanding $\gal$ action on rational functions having $\{0,1,\infty\}$ as their ramification values. This case is rich enough in that the action is faithful in genus $0$, \cite{sch-dd-on-rsphere}. Hence explicit computation of $\phi$ of Belyi pairs $(\PP^{1},\phi)$ is an indispensable part of the theory. 

A map $\psi \colon \CC \lra \CC$ of the form $\psi(z) = \lambda z + \mu$, where $\lambda,\mu \in \CC$ with $\lambda \neq 0$, is called an \emph{affine map} of $\CC$. We define a map $\phi \colon \PP^{1} \lra \PP^{1}$ to be a quotient of an affine map if there is a commutative diagram:

\begin{figure}[!h]
	\centering
	\begin{tikzpicture}[description/.style={fill=white,inner sep=1.5pt}]
		\matrix (m) [matrix of math nodes, row sep=3em,
		column sep=1.5em, text height=1.5ex, text depth=0.5ex]
		{  \CC/\Omega				&& \CC/\Omega \\ 
			\PP^{1}	&& \PP^{1} \\};
		\path[->,font=\scriptsize]
		(m-1-3) edge node [auto] {$\pi$} (m-2-3);
		\path[->,font=\scriptsize]
		(m-1-1) edge node [auto] {$\pi$} (m-2-1);
		\path[->,font=\scriptsize]
		(m-1-1) edge node [auto] {$ \psi$} (m-1-3);
		\path[->,font=\scriptsize]
		(m-2-1) edge node [auto] {$\phi$} (m-2-3);
	\end{tikzpicture}
	\label{diag:defining/X}
\end{figure}
\newpage
\noindent where $\Omega$ is a non-trivial discrete subgroup of $\CC$ and $\pi$ is a non-constant finite map. If $\Omega$ is generated by only one element, then by a change of variable, we may assume that the generator is $2\pi$. Take $\psi(z) = \lambda \,z$ where $\lambda \in \ZZ$. The commutativity of the above diagram implies that $\phi(\pi(z)) =  \pi(\lambda \, z)$. For $\pi(z) = e^{iz}$ we have $\phi(w) = w^{\lambda}$. If $\pi(z) = e^{iz} + e^{-iz}$ then $\phi(w + w^{-1}) = w^{\lambda} + w^{-\lambda}$. Functions satisfying this functional equation are called Chebyshev polynomials. For $\lambda = 2,3$, and $4$ we obtain $w^{2}-2$, $w^{3} - 3w$ and $w^{4} - 4w^{2} + 2$ respectively. In fact, for such additive subgroups, it turns out that these are the only possibilities up to conjugation, \cite[Lemma 3.8]{milnor/lattes}. We would like to point out that both $w^{\lambda}$ and Chebyshev polynomials are Belyi morphisms, i.e. ramified at most over $3$ points. 

In this work, our aim is to determine all {\lattes} maps which induce Belyi morphisms. We will show than whenever $\CC/\Omega$ does not have complex multiplication, the only isogeny that induces a Belyi morphism is multiplication by $2$ (and its translates). This gives a family of Belyi maps parametrized over the moduli space of elliptic curves,  In each of the non-CM cases, we obtain  families of Belyi morphisms, parametrized over the integers.

\section{{\lattes} maps}
\label{sec:lattes/maps}

From now on we assume that $\Omega$ is of rank two (i.e. a lattice in $\CC$), and fix the generators to be $\omega_{1}$ and $\omega_{2}$ so that $\CC/\Omega$ is an elliptic curve. In this case, the map $\phi$ is called a \emph{{\lattes} map}, \cite{lattes}. 

Since any morphism between elliptic curves is a composition of an isogeny\footnote{Recall that an isogeny between two elliptic curves is a morphism which is also a group homomorphism.} with a translation, we will limit our consideration to isogenies. For any integer $n$, multiplication by $n$ map from $\CC/\Omega$ to $\CC/\Omega$, defined as mapping each $z$ to $n \, z$\footnote{As usual, we define $[-n](z) = - [n](z)$.} is an isogeny of $\CC/\Omega$. 

By $e_{f}(x)$ let us denote the ramification index of a map $f \colon X \lra Y$ at the point $x \in X$, so that $1 \leq e_{f}(x)\leq \deg(f)$. In particular, we say that $f$ is \emph{ramified} at $x$ (or $x$ is a ramification point of $f$) if $1 < e_{f}(x)$ and $f$ is \emph{unramified} if $e_{f}(x) = 1$. If $x$ is a ramification point then $f(x)$ will be called a \emph{ramification value} of $f$. By Riemann-Hurwitz formula, morphisms from $\CC/\Omega$ to $\CC/\Omega$ are unramified, and hence $[n]$ is unramified, that is $e_{\psi}(z) = 1$ for any $z \in \CC/\Omega$. We have $e_{\psi}(z) \, e_{\pi}(\psi(z)) = e_{\pi}(z) \, e_{f} (\pi(z))$, therefore we obtain:

\begin{lemma}
	The ramification values of $f\colon \PP^{1} \lra \PP^{1}$ are a subset of the ramification values of $\pi \colon \CC/\Omega \lra \PP^{1}$.
	\label{thm:ramification/of/lattes/maps}
\end{lemma}

The following theorem limits the number of possibilities for the map \linebreak $\pi \colon \CC/\Omega \lra \PP^{1}$:
\begin{theorem}{\cite[Theorem 3.1]{milnor/lattes}}
	Let $\phi$ be a {\lattes} map. Then there is a non-trivial subgroup, denoted $\Gamma$, of automorphisms of $\CC/\Omega$ so that $\pi$ is the composition of the following maps:
	$$\CC/\Omega \lra (\CC/\Omega)/\Gamma \lra \PP^{1}.$$
\end{theorem}

Given $\Omega$, the full automorphism group of $\CC/\Omega$ is either $\ZZ/6\ZZ$ (if $j(E) = 1728$), $\ZZ/4\ZZ$ (if $j(E) = 0$) or $\ZZ/2\ZZ$ (otherwise). In the first two cases we say that $\CC/\Omega$ has complex multiplication (CM for short). One can determine the map $\pi$ in each case:

\begin{theorem}{\cite[Proposition~6.37]{silverman/dynamical/systems}}
	Given $\Omega$, the map $\pi \colon \CC/\Omega \lra \PP^{1}$ is
	$$ \pi(z) =  \begin{cases} 
								\wp(z) & \mbox{, if }\Gamma = \ZZ/2\ZZ \\
								\wp'(z) & \mbox{, if }\Gamma = \ZZ/3\ZZ \\
								(\wp(z))^{2} & \mbox{, if }\Gamma = \ZZ/4\ZZ \\
								(\wp(z))^{3} & \mbox{, if }\Gamma = \ZZ/6\ZZ
   \end{cases}
	$$
	where $\wp$, and $\wp'$ denote the Weierstra{\ss}' elliptic functions corresponding to $\Omega$.
	\label{thm:maps/to/P1}
\end{theorem}

We refer to \cite{hurwitz/allgemeine/funktionentheorie} for the following results concerning $\wp$ and $\wp'$:

\begin{theorem}
	Let $\Omega$ be a lattice and let $\wp$ be the corresponding Weierstra{\ss}' elliptic function. Then
	\begin{itemize}
		\item [i.] the equation $\wp(z_{1}) = \wp(z_{2})$ holds if and only if either $$z_{1} + z_{2} \equiv 0 \mod\Omega \mbox{ or } z_{1} - z_{2} \equiv 0 \mod \Omega,$$
		\item [ii.] the only solution of the equation $\wp'(z) = 0$ are $\frac{1}{2}\omega_{1}$, $\frac{1}{2}\omega_{2}$ and $\frac{1}{2} \omega_{3}:= \frac{1}{2}(\omega_{1} + \omega_{2})$, i.e. half-periods,
		\item [iii.] for any distinct $z_{1}, z_{2} \in \CC \sm \Omega$, we have 
			$$\wp(z_{1}+z_{2}) +\wp(z_{1}) + \wp(z_{2})  = \frac{1}{4}\left(\frac{\wp'(z_{1}) - \wp'(z_{2})}{\wp(z_{1}) - \wp(z_{2})}\right)^{2}$$ 
			\noindent whenever $z_{1}+z_{2} \notin \Omega$
	\end{itemize}
	\label{thm:facts/about/elliptic/functions}
\end{theorem}

\noindent In (iii.) taking limit as $z_{1} \lra z_{2}$ we obtain the duplication formula:
\begin{eqnarray*}
	\wp(2z) = -2\wp(z) + \frac{1}{4} \left( \frac{\wp''(z)}{\wp'(z)}   \right)
\end{eqnarray*}

If we use $E:y^{2} = x^{3} + ax + b$ model for $\CC/\Omega$, we obtain:

\begin{eqnarray}
	\wp(2z) = \frac{\wp(z)^{4} - 2z\wp(z)^{2} - 8b\wp(z) + a^{2}}{4(\wp(z)^{3} + z\wp(z) + b)}
	\label{eq:duplication/formula}
\end{eqnarray}

\noindent which is valid whenever $z$ and $2z$ are not in $\Omega$.

\subsection{Non-CM case.}
\label{sec:non/CM}

Assuming $\CC/\Omega$ does not have CM, the only possible non-trivial automorphism group is $\ZZ/2\ZZ$. In this case, the map $\pi$ is $\wp$, and hence is of degree $2$. So the candidates for ramification values of $\phi$ are $\infty$, and values of $\wp$ at half-periods. As discussed above, we take $\psi$ to be $[n]$. 

Say $n$ is an odd integer. Each half period $\frac{\omega_{i}}{2} \in \CC/\Omega$ and $0$ has $n^{2}$ pre-images in $\CC/\Omega$. Among these pre-images, $n^{2}-1$ of them come in pairs so that the value of $\wp$ at those pairs agree. By Riemann-Hurwitz formula we see that this is all the ramification:
$$-2 = n^{2}(-2) + 4 (\frac{n^{2}-1}{2})$$
\noindent So we conclude that the map $\phi$ is not a Belyi morphism.

Suppose now that $n$ is an even integer. For each of the half periods, to every point $z \in [n]^{-1}(\frac{\omega_{i}}{2})$ there corresponds an element, say $z_{o} \in [n]^{-1}(\frac{\omega_{i}}{2})$ so that $z + z_{o} \equiv 0 \mod \Omega$. For $0$, among the set $[n]^{-1}(0)$ apart from $o,\omega_{1}/2,\omega_{2}/2$ and $\omega_{3}/2$ to every point in $z \in [n]^{-1}(0)$ there is an element $z_{o}$ so that $z + z_{o} \equiv 0 \mod \Omega$. Riemann-Hurwitz formula reads:
$$-2 = n^{2}(-2) + 3 \frac{n^{2}}{2} + \frac{n^{2} - 4}{2}.$$

For the remaining case, note that $[2]^{-1}(0) = \{0,\omega_{1}/2,\omega_{2}/2,\omega_{3}/\}$. Therefore the map $\phi$ is not ramified at $\wp(0) = \infty$, i.e. ramified over $3$ points. Equivalently, $0$  is not a ramification point as $\frac{n^{2} - 4}{2} = 0$ for $n=2$. So we proved:

\begin{theorem}
	If $\CC/\Omega$ does not have CM, then the corresponding {\lattes} maps are Belyi morphisms if and only if the isogeny is multiplication by $2$ map(, or any of its translate).
	\label{thm:nonCM/case}
\end{theorem}

We would like to remark that as far as we know this is the first instance where one has a continuous family of Belyi morphisms, parametrized by $a$ and $b$ as shown in Equation~\ref{eq:duplication/formula}.

The following result can also be deduced from the proof of the previous theorem:

\begin{lemma}
	Given the commutative diagram:
	\begin{figure}[!h]
	\centering
	\begin{tikzpicture}[description/.style={fill=white,inner sep=1.5pt}]
		\matrix (m) [matrix of math nodes, row sep=3em,
		column sep=1.5em, text height=1.5ex, text depth=0.5ex]
		{  \CC/\Omega				&& \CC/\Omega \\ 
			\PP^{1}	&& \PP^{1} \\};
		\path[->,font=\scriptsize]
		(m-1-3) edge node [auto] {$\pi$} (m-2-3);
		\path[->,font=\scriptsize]
		(m-1-1) edge node [auto] {$\pi$} (m-2-1);
		\path[->,font=\scriptsize]
		(m-1-1) edge node [auto] {$[n]$} (m-1-3);
		\path[->,font=\scriptsize]
		(m-2-1) edge node [auto] {$\phi$} (m-2-3);
	\end{tikzpicture}
\end{figure}
	
	\noindent where the map $\pi$ is $(\wp(z))^{j}$ with $j = 1,2,3$\footnote{We do not lose any generality in making such a choice for $j$.}, then $\infty$ is a ramification value of $\varphi$ if and only if $n>2$ or $n<-2$.
	\label{thm:ramification/for/wp}
\end{lemma}

\subsection{CM case.}

In the CM case there are 3 possibilities for the non-trivial subgroup $\Gamma$ of the automorphism group $\Gamma$, namely, $\ZZ/3\ZZ$, $\ZZ/4\ZZ$ or $\ZZ/6\ZZ$. In what follows, we fix the zeroes of $\wp(z)$ as $d_{1}$ and $d_{2}$ in $\CC/\Omega$ so that $d_{1} + d_{2} \equiv 0 \mod\Omega$\footnote{Interested reader may consult \cite{eichler/zagier/zeroes/of/wp} for a formula for the zeros in terms of a modular type integral, and \cite{duke/imamoglu/zeros/of/wp} for a formula in terms of the modular variable.}.

\subsubsection{$\Gamma = \ZZ/3\ZZ$.} 
\begin{figure}[h!]
	\centering
	\includegraphics{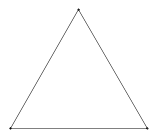}
	 \qquad
	\includegraphics{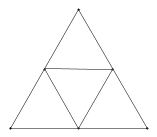}
	\qquad
	\includegraphics{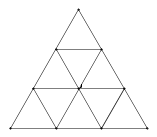}
	\caption{First three subdivisions of an equilateral triangle.}
	\label{fig:triangle/and/subdivisions}
\end{figure}

We have $\pi = \wp'$ in this case. On $\CC/\Omega$, we have $\wp''(z) = 6(\wp(z))^{2}$ hence there are at most three ramification points of the function of $\pi$, namely $0,d_{1}$ and $d_{2}$. Therefore, for any integer $n$, the map $\phi$ is ramified at most over 3 points, and hence a Belyi morphism. In fact, for any integer $n$, the \emph{dessin} corresponding to these Belyi maps can be obtained by gluing two copies of an equilateral triangle whose edges are divided into $n$ segments along their boundaries, see Figure~\ref{fig:triangle/and/subdivisions}. 

The corresponding Belyi maps can also be computed by an iteration of the duplication formula. In particular, under suitable normalization, they are given as:

\begin{eqnarray*}
	\frac{1}{8} \, \frac{(z-1)(z+1)^{3}}{(z+\frac{1}{2})^{3}} & \mbox{, whenever $n = 2$} \\
	\frac{1}{27} \frac{(z^{3}+3z^{2}-6z+1)^{3}}{z (z-1)(z^{2}-z+1)^{3}} & \mbox{, whenever $n = 3$} \\ 
	\frac{1}{8}\frac{z(z^{5}+8z^{4}-32z^3+28z^2-10z+3)^{3}}{(2z^{5} - 5z^{4}+15z^{3}-16z^{2}+4z+1/2)}& \mbox{, whenever $n = 4$}
\end{eqnarray*}

\subsubsection{$\Gamma = \ZZ/4\ZZ$.} 
As a result of Theorem~\ref{thm:maps/to/P1}, the map $\pi$ is $(\wp(z))^{2}$. The other ramification points of $\phi$ may come from the roots of $\wp$ and roots of $\wp'$, namely, zeroes of $\wp'$ which are half periods, $\frac{\omega_i}{2}$, $i = 1,2,3$ and zeroes of $\wp$ which we denote by $d_{1}$ and $d_{2}$. In this case, however, we have $d_{1} = d_{2} = e_{3}$ and $\wp(e_{1}) = -\wp(e_{2}$. These equalities may be deduced from the extra symmetry $\wp(\sqrt{-1} \, z = -\wp(z)$. The following table summarizes ramification data:

\begin{table}[h!]
	\centering
	\begin{tabular}{|c|c|c|c|}
		\hline
		parity of $n$ 	& $|\phi^{-1}(0)|$ 			& $|\phi^{-1}(\wp(e_{1}))^{2}|$ & $|\phi^{-1}(\infty)|$ \\\hline
		 even 			& $\frac{n^{2}}{4}$ 			& $\frac{2n^{2}}{4}$ & $ \frac{n^{2}-4}{4} + 3$\\\hline
		 odd 				& $\frac{n^{2}-1}{4} + 1$	& $\frac{2n^{2}-2}{4}+1$ & $\frac{n^{2}-1}{4} + 1$ \\\hline
	\end{tabular}
\end{table}

\subsubsection{$\Gamma = \ZZ/6\ZZ$.} In this case, $\pi  = \wp^{3}$. As in the previous case, the ramification points of $\pi$ come from the roots of $\wp$ and $\wp'$. This time, $(\wp(\omega_{i}/2)^{3} = (\wp(\omega_{j}/2)^{3}$ for all $i,j \in \{1,2,3\}$. Therefore, for any integer $n \in \ZZ$, the corresponding {\lattes} map $\phi\colon\PP^{1} \lra \PP^{1}$ is ramified at most over $0,\infty$ and $(\wp(\omega_{1}/2))^{3}$.

The following theorem summarizes all the results above:

\begin{theorem}
	If $\phi \colon \PP^{1} \lra \PP^{1}$ is a {\lattes} map obtained from an elliptic curve, $\CC/\Omega$, having CM by using multiplication by $n$ map on the elliptic curve, then $\phi$ is a Belyi morphism.
\end{theorem}

\bibliographystyle{abbrv}

\end{document}